\documentclass[12pt]{article}
\usepackage{amssymb, amsmath, latexsym, a4}
\usepackage[latin1]{inputenc}
\usepackage[all]{xy}

\newcommand{\rdg}{\hfill $\Box $}

\newtheorem{De}{Definition}[section]
\newtheorem{Th}[De]{Theorem}
\newtheorem{Pro}[De]{Proposition}
\newtheorem{Le}[De]{Lemma}
\newtheorem{Co}[De]{Corollary}
\newtheorem{Rem}[De]{Remark}
\newtheorem{Ex}[De]{Example}

\begin{document}

\centerline{\Large On universal central extensions of Hom-Leibniz algebras}
\bigskip \bigskip

\centerline{José Manuel Casas$^{(a)}$, Manuel Avelino Insua$^{(a)}$  and Natalia Pacheco Rego$^{(b)}$}
\bigskip \bigskip

{ \centerline{$^{(a)}$ Dpto.  Matemática Aplicada I, Univ. de Vigo, E. E. Forestal,}

\centerline{ 36005 Pontevedra, Spain}

\centerline{e-mail addresses: {\tt jmcasas@uvigo.es}, {\tt ainsua@uvigo.es}}
\medskip

 \centerline{$^{(b)}$ IPCA, Dpto. de Ciências, Campus do IPCA,
 Lugar do Aldão}
\centerline{4750-810 Vila Frescainha, S. Martinho, Barcelos,
 Portugal}}

 \centerline{e-mail address: {\tt nrego@ipca.pt} }
\bigskip

{\bf Abstract:}
In the category of Hom-Leibniz algebras we introduce the notion of representation as  adequate coefficients to construct the chain complex  to compute the Leibniz homology of Hom-Leibniz algebras. We study universal central extensions of Hom-Leibinz algebras and generalize some classical results, nevertheless it is necessary to introduce  new notions of $\alpha$-central extension, universal $\alpha$-central extension and $\alpha$-perfect Hom-Leibniz algebra. We prove that an $\alpha$-perfect Hom-Lie algebra  admits a universal $\alpha$-central extension in the categories of Hom-Lie and Hom-Leibniz algebras and we obtain the relationships between both.
In case  $\alpha = Id$ we recover the corresponding results on universal central extensions of Leibniz algebras.
\bigskip

\bigskip


\section{Introduction}

The Hom-Lie algebra structure was initially introduced in \cite{HLS}  motivated by examples of deformed Lie algebras coming from twisted discretizations of vector fields. Hom-Lie algebras are $\mathbb{K}$-vector spaces endowed with a bilinear skew-symmetric bracket satisfying a  Jacobi identity twisted by a map. When this map is the identity map, then the definition of Lie algebra is recovered.

From the introductory paper  \cite{HLS}, this algebraic structure and other related ones as Hom-associative, Hom-Leibniz and Hom-Nambu algebras were studied in several papers \cite{Is, MS, MS2, MS3, Sh, Yau 1, Yau, Yau 2} and references given therein.

Following the generalization in \cite{Lo} from Lie to Leibniz algebras, it is natural to describe same generalization in the framework of Hom-Lie algebras. In this way, the notion of Hom-Leibniz algebra was firstly introduced in \cite{MS} as $\mathbb{K}$-vector spaces $L$ together with a linear map $\alpha: L \to L$, endowed with a bilinear bracket operation $[-,-] : L\times L \to L$ which satisfies the Hom-Leibniz identity $[\alpha(x),[y,z]]=[[x,y],\alpha(z)]-[[x,z],\alpha(y)]$, for all $x, y, z \in L$, and it was the subject of the recent papers \cite{AMM, Is, MS2, MS3, Yau}. A (co)homology theory and an initial study of universal central extensions was given in \cite{ChS}.

Our goal in the present paper is to generalize properties and characterizations of universal central extensions of Leibniz algebras in \cite{CC, CV, Gn 1} to Hom-Leibniz algebras setting. But an important fact, which is  the composition of two central extension is not central as the counterexample \ref{counterexample} below shows, doesn't allow us to obtain a complete generalization of the classical results. Nevertheless it leads us to introduce new notions as $\alpha$-centrality or $\alpha$-perfection in order to generalize classical results. On the other hand, we prove that an $\alpha$-perfect Hom-Lie algebra admits a universal $\alpha$-central extension in the categories of Hom-Lie and Hom-Leibnz algebras, then  one of our main results establishes the relationships between both universal $\alpha$-central extensions. When we rewrite these relationships for Lie and Leibniz algebras, i.e. the twisting endomorphism is the identity, we recover the corresponding results given in \cite{Gn 1}.

In order to achieve our goal we organize the paper as follows: section 2 is dedicated to introduce the background material on Hom-Leibniz algebras. In section 3 we introduce  co-representations, which are the adequate coefficients to define the chain complex from which we compute the homology of a Hom-Leibniz algebra with coefficients. Low-dimensional homology $\mathbb{K}$-vector spaces are obtained. In case $\alpha = Id$   we recover the homology of Leibniz algebras  \cite{Lo, LP}.
 In section 4 we present our main results on universal central extensions, namely we extend classical results and present a counterexample showing that the composition of two central extension is not a central extension. This fact lead us to define $\alpha$-central extensions: an extension $\pi : (K,\alpha_k) \twoheadrightarrow (L,\alpha_L)$ is said to be $\alpha$-central if  $[\alpha({\rm Ker}\ (\pi)), K] = 0 = [K,\alpha({\rm Ker}\ (\pi))] $, and is said to be central if  $[{\rm Ker}\ (\pi), K] = 0 = [K,{\rm Ker}\ (\pi)]$. Clearly central extension implies $\alpha$-central extension and both notions coincide in case $\alpha = Id$.

 We can extend classical results on universal central extensions of Leibniz algebras in \cite{CC, CV, Gn 1} to the Hom-Leibniz algebras setting as: a Hom-Lie algebra is perfect if and only if admits a universal central extension and the kernel of the universal central extension is the second homology with trivial coefficients of the Hom-Leibniz algebra. Nevertheless, other result as: if a central extension $0 \to (M, \alpha_M) \stackrel{i} \to (K,\alpha_K) \stackrel{\pi} \to (L, \alpha_L) \to 0$ is universal, then $(K,\alpha_K)$ is perfect and every central extension of $(K,\alpha_K)$ is split only holds for universal $\alpha$-central extensions, which means that only lifts on $\alpha$-central extensions. Other relevant result, which cannot be extended in the usual way, is: if $0 \to (M, \alpha_M) \stackrel{i} \to (K,\alpha_K) \stackrel{\pi} \to (L, \alpha_L) \to 0$ is a universal   $\alpha$-central extension, then $H_1^{\alpha}(K) = H_2^{\alpha}(K) = 0$. Of course, when the twisting endomorphism is the identity morphism, then all the new notions and all the new results coincide with the classical ones.

In section 5 we prove that an $\alpha$-perfect Hom-Lie algebra $(L,[-,-],\alpha_L)$, that is $L=[\alpha_L(L),\alpha_L(L)]$, admits a universal $\alpha$-central extension in the categories of Hom-Lie and Hom-Leibniz algebras, then we obtain the relationships between both. Our main results in this section generalize the relationships between the universal central extensions of a Lie algebra in the categories of Lie and Leibniz algebras given in \cite{Gn} when we consider Leibniz algebras as Hom-Leibniz algebras, i.e. when the twisting endomorphism is the identity.

\section{Hom-Leibniz algebras}
In this section we introduce necessary material on  Hom-Leibniz algebras which will be used in subsequent sections.

\begin{De}\label{HomLeib} \cite{MS}
A Hom-Leibniz algebra is a triple $(L,[-,-],\alpha_L)$ consisting of a $\mathbb{K}$-vector space $L$, a bilinear map $[-,-] : L \times L \to L$ and a $\mathbb{K}$-linear map $\alpha_L : L \to L$ satisfying:
\begin{equation} \label{def}
 [\alpha_L(x),[y,z]]=[[x,y],\alpha_L(z)]-[[x,z],\alpha_L(y)] \ \ ({\rm Hom-Leibniz\ identity})
\end{equation}
for all $x, y, z \in L$.
\end{De}

In terms of the adjoint representation $ad_x : L  \to L, ad_x(y)=[y,x]$, the Hom-Leibniz identity can be written as follows: $$ad_{\alpha(z)} \cdot ad_y = ad_{\alpha(y)} \cdot ad_z + ad_{[y,z]} \cdot \alpha$$

\begin{De} \cite{Yau 1}
A Hom-Leibniz algebra $(L,[-,-],\alpha_L)$  is said to be multiplicative  if the $\mathbb{K}$-linear map $\alpha_L$ preserves the bracket, that is, if $\alpha_L [x,y] = [\alpha_L(x),\alpha_L(y)]$, for all $x, y \in L$.
\end{De}

\begin{Ex}\label{ejemplo 1} \
\begin{enumerate}
\item[a)] Taking $\alpha = Id$ in Definition \ref{HomLeib} we obtain the definition of  Leibniz algebra \cite{Lo}. Hence Hom-Leibniz algebras include Leibniz algebras as a full subcategory, thereby motivating the name "Hom-Leibniz algebras" as a deformation of Leibniz algebras twisted by a homomorphism. Moreover it  is a multiplicative Hom-Leibniz algebra.

    \item[b)] Hom-Lie algebras \cite{HLS} are Hom-Leibniz algebras whose bracket satisfies the condition $[x,x]=0$, for all $x$. So Hom-Lie algebras can be considered as a full subcategory of Hom-Leibniz algebras category. For any multiplicative Hom-Leibniz algebra $(L,[-,-],\alpha_L)$ there is associated the Hom-Lie algebra $(L_{\rm Lie},[-,-],\widetilde{\alpha})$, where $L_{\rm Lie} = L/L^{\rm ann}$, the bracket is the canonical bracket induced on the quotient and  $\widetilde{\alpha}$ is the homomorphism naturally induced by $\alpha$. Here $L^{\rm ann} = \langle \{[x,x] : x \in L \} \rangle$.

    \item[b)] Let $(D,\dashv,\vdash, \alpha_D)$ be a Hom-dialgebra. Then $(D,\dashv,\vdash, \alpha_D)$ is a  Hom-Leibniz algebra with respect to the bracket  $[x,y]=x \dashv y - y \vdash x$, for all $x,y \in A$ \cite{Yau}.

  \item[c)] Let $(L,[-,-])$ be a Leibniz algebra and $\alpha_L:L \to L$  a Leibniz algebra endomorphism. Define $[-,-]_{\alpha} : L \otimes L \to L$ by $[x,y]_{\alpha} = [\alpha(x),\alpha(y)]$, for all $x, y \in L$. Then $(L,[-,-]_{\alpha}, \alpha_L)$ is a multiplicative Hom-Leibniz algebra.

\item[d)] Abelian or commutative Hom-Leibniz algebras are $\mathbb{K}$-vector spaces $L$ with trivial bracket and any linear map $\alpha_L :L \to L$.

\end{enumerate}
\end{Ex}

\begin{De}\label{homo}
A homomorphism of  Hom-Leibniz algebras $f:(L,[-,-],\alpha_L) \to (L',[-,-]',\alpha_{L'})$ is a $\mathbb{K}$-linear map $f : L \to L'$ such that
\begin{enumerate}
\item[a)] $f([x,y]) =[f(x),f(y)]'$
\item [b)] $f \cdot \alpha_L(x) = \alpha_{L'} \cdot  f(x)$
\end{enumerate}
for all $x, y \in L$.
\end{De}

A homomorphism of multiplicative Hom-Leibniz algebras is a homomorphism of the underlying Hom-Leibniz algebras.
\bigskip

So we have defined the category ${\sf Hom-Leib}$ (respectively, ${\sf Hom-Leib_{\rm mult}})$ whose objects are Hom-Leibniz (respectively, multiplicative Hom-Leibniz) algebras and whose morphisms are the homomorphisms of Hom-Leibniz (respectively, multiplicative Hom-Leibniz) algebras.
There is an obvious inclusion functor $inc : {\sf Hom-Leib_{\rm mult}} \to {\sf Hom-Leib}$. This functor has as left adjoint the multiplicative functor $(-)_{\rm mult} : {\sf Hom-Leib} \to {\sf Hom-Leib_{\rm mult}}$ which assigns to a Hom-Leibniz algebra $(L,[-,-],\alpha_L)$ the multiplicative Hom-Leibniz algebra $(L/I,[-,-],\tilde{\alpha})$, where $I$ is the two-sided ideal of $L$ spanned by the elements $\alpha_L[x,y]-[\alpha_L(x),\alpha_L(y)]$, for all $x, y \in L$.
\bigskip

In the sequel we refer Hom-Leibniz algebra to a multiplicative Hom-Leibniz algebra and we shall use the shortened notation $(L,\alpha_L)$ when there is not confusion with the bracket operation.

\begin{De}
Let $(L,[-,-],\alpha_L)$ be a Hom-Leibniz algebra. A  Hom-Leibniz subalgebra $H$ is a linear subspace of $L$, which is closed for the bracket and invariant by $\alpha_L$, that is,
\begin{enumerate}
\item [a)] $[x,y] \in H,$ for all $x, y \in H$
\item [b)] $\alpha_L(x) \in H$, for all $x \in H$
\end{enumerate}

A  Hom-Leibniz subalgebra $H$ of $L$ is said to be a  two-sided Hom-ideal if $[x,y], [y,x] \in H$, for all $x \in H, y \in L$.

If $H$ is a two-sided  Hom-ideal of $L$, then the quotient $L/H$ naturally inherits a structure of Hom-Leibniz algebra, which is said to be the quotient Hom-Leibniz algebra.
\end{De}

\begin{De}
Let $H$ and $K$ be two-sided Hom-ideals of a Hom-Leibniz algebra $(L,[-,-],\alpha_L)$. The commutator  of $H$ and $K$, denoted by $[H,K]$, is the Hom-Leibniz subalgebra of $L$ spanned by the brackets $[h,k], h \in H, k \in K$.
\end{De}

Obviously, $[H,K] \subseteq H \cap K$ and $[K,H] \subseteq H \cap K$. When $H = K =L$, we obtain the definition of derived Hom-Leibniz subalgebra. Let us observe that, in general, $[H,K]$ is not   a Hom-ideal, but if $H, K \subseteq \alpha_L(L)$, then $[H,K]$ is a two-sided ideal of $\alpha_L(L)$. When $\alpha = Id$, the classical notions are recovered.

\begin{De}
Let $(L,[-,-],\alpha_L)$ be a Hom-Leibnz algebra. The subspace $Z(L) = \{ x \in L \mid [x, y] =0 = [y,x], \text{for\ all}\ y \in L \}$ is said to be the center of $(L,[-,-],\alpha_L)$.

When $\alpha_L : L \to L$ is a surjective homomorphism, then  $Z(L)$ is a Hom-ideal of $L$.
\end{De}

\section{Homology of Hom-Leibniz algebras}

In this section, we introduce the notion of Hom-co-representation, construct a chain complex from which we define the homology $\mathbb{K}$-vector spaces of a Hom-Leibniz algebra with coefficients on a Hom-co-representation and we  interpret low-dimensional homology  $\mathbb{K}$-vector spaces as well.

\begin{De}
Let $(L,[-,-],\alpha_L)$ be a Hom-Leibniz algebra. A Hom-co-represen\-tation of $(L,[-,-],\alpha_L)$ is a $\mathbb{K}$-vector space $M$ together with two bilinear maps $\lambda : L \otimes M \to M, \lambda(l\otimes m) =l \centerdot m$, and $\rho : M \otimes L \to M, \rho(m \otimes l) = m \centerdot l$, and a $\mathbb{K}$-linear map $\alpha_M : M \to M$ satisfying the following identities:
\begin{enumerate}
\item $[x,y] \centerdot \alpha_M(m) = \alpha_L(x) \centerdot (y \centerdot m) - \alpha_L(y) \centerdot (x \centerdot m)$.

\item $\alpha_L(y) \centerdot (m \centerdot x) = (y \centerdot m) \centerdot \alpha_L(x) - \alpha_M(m) \centerdot [x,y]$.

\item $(m \centerdot x) \centerdot \alpha_L(y) = \alpha_M(m) \centerdot [x,y] - (y \centerdot m) \centerdot \alpha_L(x)$.

\item $\alpha_M(x \centerdot m) = \alpha_L(x) \centerdot \alpha_M(m) $

\item $\alpha_M(m \centerdot x) =  \alpha_M(m) \centerdot \alpha_L(x)$
\end{enumerate}
for any $x, y \in L$ and $m \in M$
\end{De}

From the second and third identities above directly follows
 \begin{equation}
 \alpha_L(y) \centerdot (m \centerdot x) + (m \centerdot x) \centerdot \alpha_L(y)=0
 \end{equation}

\begin{Ex}\ \label{ejemplo}
\begin{enumerate}
\item[a)] Let $M$ be a  co-representation of a Leibniz algebra L \cite{LP}. Then $(M,Id_M)$ is a  Hom-co-representation of the  Hom-Leibniz algebra $(L,Id_L)$.

\item[b)] Every  Hom-Leibniz algebra  $(L,[-,-],\alpha_L)$ has a Hom-co-representation structure on itself given by the actions $$x \centerdot m =-[m, x]; \quad m \centerdot x = [m, x]$$
    where $x \in L$ and $m$ is an element of the underlying $\mathbb{K}$-vector space to $L$.
\end{enumerate}
\end{Ex}

 Let $(L,[-,-],\alpha_L)$ be a Hom-Leibniz algebra and $(M,\alpha_M)$ be a Hom-co-representa\-tion of $(L,[-,-],\alpha_L)$. Denote $CL_n^{\alpha}(L,M) :=M \otimes L^{\otimes n}, n \geq 0$. We define the $\mathbb{K}$-linear map $$d_n: CL_n^{\alpha}(L,M) \to CL_{n-1}^{\alpha}(L,M)$$ by $$d_n(m \otimes x_1 \otimes \dots \otimes x_n)= m \centerdot x_1 \otimes \alpha_L(x_2) \otimes \dots \otimes \alpha_L(x_n) +$$ $$\sum_{i=2}^n (-1)^i x_i \centerdot m \otimes \alpha_L(x_1) \otimes \dots \otimes \widehat{\alpha_L(x_i)} \otimes \dots \otimes \alpha_L(x_n) +$$ $$\sum _{1 \leq i < j \leq n} (-1)^{j+1} \alpha_M(m) \otimes \alpha_L(x_1) \otimes \dots \otimes \alpha_L(x_{i-1})\otimes [x_i,x_j] \otimes \dots \otimes \widehat{\alpha_L(x_j)} \otimes \dots \otimes \alpha_L(x_n)$$

The chain complex $(CL_{\star}^{\alpha}(L,M), d_{\star})$ is  well-defined, that is, $d_n \cdot d_{n+1} =0, n \geq 0$. Indeed,
if we define for any $y\in L$ and $n\in \mathbb{N}$ two $\mathbb{K}$-linear maps,
 $$\theta_{n}\left(  y\right)  : CL_{n}^{\alpha
}\left(  L,M\right)  \longrightarrow CL_{n}^{\alpha}\left(  L,M\right) $$ given by
$$\theta_{n}\left(  y\right)   \left(
m\otimes x_{1}\otimes\cdots\otimes x_{n}\right)  =-y\cdot m\otimes\alpha
_{L}\left(  x_{1}\right)  \otimes\cdots\otimes\alpha_{L}\left(  x_{n}\right) +
$$
$$\overset{n}{\underset{i=1}{\displaystyle\sum}} \alpha_{M}\left(  m\right)  \otimes\alpha_{L}\left(  x_{1}\right)
\otimes\cdots\otimes\left[  x_{i},y\right]  \otimes\cdots\otimes\alpha
_{L}\left(  x_{n}\right)$$
and
$$i_{n}\left( \alpha_L( y) \right)  : CL_{n}^{\alpha}\left(  L,M\right)  \longrightarrow
CL_{n+1}^{\alpha}\left(  L,M\right)$$ given by
$$i_{n}\left(  \alpha_{L}\left(
y\right)  \right)  \left(  m\otimes x_{1}\otimes\cdots\otimes x_{n}\right)
=\left(  -1\right)  ^{n}m\otimes x_{1}\otimes\cdots\otimes x_{n}\otimes y$$
then the following formulas hold:

\begin{Pro}(Generalized  Cartan's Formulas) \label{CartanLeib}

The following identities hold:
\begin{enumerate}
\item[a)] $d_{n+1} \cdot i_{n}\left(  \alpha_{L}\left(  y\right)  \right)
+i_{n-1} \left(  \alpha_{L}^{2}\left(  y\right)  \right) \cdot  d_{n}
=\theta_{n}\left(  y\right),$ for  $n \geq 1.$

\item[b)] $\theta_{n}\left(  \alpha_{L}\left(  x\right)  \right) \cdot
\theta_{n}\left(  y\right)  -\theta_{n} \left(  \alpha_{L}\left(  y\right)
\right) \cdot  \theta_{n}\left(  x\right)  =-\theta_{n} \left(  \left[
x,y\right]  \right) \cdot \left(  \alpha_{M}\otimes\alpha_{L}^{\otimes n}\right),$
for $n \geq 0.$

\item[c)] $\theta_{n}\left(  x\right) \cdot  i_{n-1}\left( \alpha_L( y) \right)
-i_{n-1}\left(  \alpha_{L}^2\left(  y\right)  \right) \cdot  \theta_{n-1}\left(
x\right)  =i_{n-1}(\alpha_L(\left[  y,x\right])) \cdot \left(  \alpha_{M}\otimes\alpha
_{L}^{\otimes n -1}\right),$ for $n \geq 1.$

\item[d)] $\theta_{n-1}\left(  \alpha_{L}\left(  y\right)  \right) \cdot
d_{n}=d_{n} \cdot \theta_{n}\left(  y\right),$ for $n \geq 1.$

\item[e)] $d_n \cdot d_{n+1}=0,$ for $n \geq 1$.
\end{enumerate}
\end{Pro}
{\it Proof.} The proof is followed by a standard mathematical induction, so we omit it. \rdg
\bigskip

 The homology of the chain complex $(CL_{\star}^{\alpha}(L,M), d_{\star})$ is called de homology of the Hom-Leibniz algebra  $(L,[-,-],\alpha_L)$ with coefficients in the Hom-co-representation $(M,\alpha_M)$ and we will denote it by:
 $$HL_{\star}^{\alpha}(L,M) :=H_{\star}(CL_{\star}^{\alpha}(L,M), d_{\star})$$

 Now we are going to compute low dimensional homologies. So, for $n=0$, a direct checking shows that $$HL^{\alpha}_0(L,M) = \frac{M}{M_L}$$
 where $M_L=\{m \centerdot l: m \in M, l \in L\}$.

 If $M$ is a trivial Hom-co-representation, that is, $m \centerdot l = l \centerdot m = 0$, then $$HL^{\alpha}_1(L,M) = \frac{M \otimes L}{\alpha_M(M) \otimes [L,L]}$$

\begin{Pro}\label{iso}
Let $(L, [-,-], \alpha_L)$ be a Hom-Leibniz algebra, which is considered as a Hom-co-representation of itself as in Example \ref{ejemplo} b) and $\mathbb{K}$ as a trivial Hom-co-representation of $(L,[-,-],\alpha_L)$. Then
$$HL_{n}^{\alpha} (L,L) \cong HL_{n + 1}^{\alpha}(L,\mathbb{K}), n \geq 0$$
\end{Pro}
{\it Proof.} Obviously $- Id: CL_n^{\alpha}(L,L) \to CL_{n+1}^{\alpha}(L,\mathbb{K}), n \geq 0,$ defines a chain isomorphism, hence the isomorphism in the homologies. \rdg

\begin{Rem}
Let $L$ be a Leibniz algebra. If we consider $L$ as a Hom-Leibniz algebra as in Example \ref{ejemplo 1} a) and we rewrite Proposition \ref{iso} with this particular case, considering $(L,\alpha_L)$ with a Hom-co-representation structure as in Example \ref{ejemplo} b), then  we recover isomorphism 6.5 in \cite{Lo} (see also application 3.1 in \cite{Gn}).
\end{Rem}

\section{Universal central extensions}

Through this section we deal with universal central extensions of Hom-Leibniz algebras. We generalize classical results on universal central extensions of Leibniz algebras, but an inconvenient, which is the composition of central extensions is not central as Example \ref{counterexample} shows,  doesn't allow us to obtain the complete generalization of all classical results. Nevertheless, this fact lead us to introduce a new concept of centrality.

\begin{De} \label{alfacentral}
A short exact sequence of Hom-Leibniz algebras $(K) : 0 \to (M, \alpha_M) \stackrel{i} \to (K,\alpha_K) \stackrel{\pi} \to (L, \alpha_L) \to 0$ is said to be central if $[M, K] = 0 = [K, M]$. Equivalently,  $M \subseteq Z(K)$.

We say that $(K)$ is $\alpha$-central if $[\alpha_M(M), K] = 0 = [K,\alpha_M(M)]$. Equivalently, $\alpha_M(M) \subseteq Z(K)$.
\end{De}

\begin{Rem}
Obviously, every central extension is an $\alpha$-central extension.
Note that in the case $\alpha_M = Id_M$,  both notions coincide.
\end{Rem}

\begin{De} \label{Def universal}
A central extension $(K) : 0 \to (M, \alpha_M) \stackrel{i} \to (K,\alpha_K) \stackrel{\pi} \to (L, \alpha_L) \to 0$ is said to be universal if for every central extension $(K') : 0 \to (M', \alpha_{M'}) \stackrel{i'} \to (K',\alpha_{K'}) \stackrel{\pi'} \to (L, \alpha_L) \to 0$ there exists a unique homomorphism of Hom-Leibniz algebras  $h : (K,\alpha_K) \to (K',\alpha_{K'})$ such that $\pi'\cdot h = \pi$.

We say that the central extension  $(K) : 0 \to (M, \alpha_M) \stackrel{i} \to (K,\alpha_K) \stackrel{\pi} \to (L, \alpha_L) \to 0$ is universal $\alpha$-central if for every $\alpha$-central extension $(K) : 0 \to (M', \alpha_{M'}) \stackrel{i'} \to (K',\alpha_{K'}) \stackrel{\pi'} \to (L, \alpha_L) \to 0$ there exists a unique homomorphism of Hom-Leibniz algebras  $h : (K,\alpha_K) \to (K',\alpha_{K'})$ such that $\pi'\cdot h = \pi$.
\end{De}

\begin{Rem} \label{rem}
Obviously, every universal $\alpha$-central extension is a universal central extension.
Note that in the case $\alpha_M = Id_M$,  both notions coincide.
\end{Rem}

\begin{De}
We say that a Hom-Leibniz algebra $(L, \alpha_L)$ is perfect if $L = [L,L]$.
\end{De}

\begin{Le} \label{lema 1}
Let $\pi : (K,\alpha_K) \to (L, \alpha_L)$ be a surjective homomorphism of Hom-Leibniz algebras. If $(K,\alpha_K)$ is a perfect Hom-Leibniz algebra, then $(L, \alpha_L)$  is a perfect Hom-Leibniz algebra as well.
\end{Le}

\begin{Le} \label{lema 2}
Let $0 \to (M, \alpha_M) \stackrel{i} \to (K,\alpha_K) \stackrel{\pi} \to (L, \alpha_L) \to 0$ be an $\alpha$-central extension and $(K,\alpha_K)$ a perfect Hom-Leibniz algebra. If there exists a homomorphism of Hom-Leibniz algebras $f : (K,\alpha_K) \to (A, \alpha_A)$ such that $\tau \cdot f = \pi$, where $0 \to (N, \alpha_N) \stackrel{j} \to (A,\alpha_A) \stackrel{\tau} \to (L, \alpha_L) \to 0$ is a central extension, then $f$ is unique.
\end{Le}

The proofs of these two last Lemmas use classical arguments, so we omit it.

\begin{Le} \label{lema 3}
If $0 \to (M, \alpha_M) \stackrel{i} \to (K,\alpha_K) \stackrel{\pi} \to (L, \alpha_L) \to 0$ is a universal central extension, then $(K,\alpha_K)$ and $(L, \alpha_L)$ are perfect Hom-Leibniz algebras.
\end{Le}
{\it Proof.} Assume that $(K,\alpha_K)$ is not a perfect Hom-Leibniz algebra, then $[K,K] \varsubsetneq K$. Consider $I$ the smallest Hom-ideal generated by $[K,K]$, then $(K/I, \widetilde{\alpha})$, where $\widetilde{\alpha}$ is the induced natural homomorphism, is an abelian Hom-Leibniz algebra and, consequently, is a trivial Hom-co-representation of $(L, \alpha_L)$. Consider the central extension $0 \to (K/I, \widetilde{\alpha}) \to (K/I \times L, \widetilde{\alpha}\times \alpha_L) \stackrel{pr}\to  (L, \alpha_L) \to 0$, then the homomorphisms of Hom-Leibniz algebras $\varphi, \psi : (K,\alpha_K) \to (K/I \times L, \widetilde{\alpha}\times \alpha_L)$ given by $\varphi(k)=(k+I,\pi(k))$ and $\psi(k)=(0,\pi(k)), k \in K$ verify that $pr\cdot  \phi = \pi = pr\cdot  \psi$, so $0 \to (M, \alpha_M) \stackrel{i} \to (K,\alpha_K) \stackrel{\pi} \to (L, \alpha_L) \to 0$ cannot be a universal central extension.

Lemma \ref{lema 1} ends the proof. \rdg
\bigskip

Classical categories as groups, Lie algebras, Leibniz algebras and other similar
ones share the following property: the composition of two central extensions is
a central extension as well. This property is absolutely necessary in order to obtain characterizations
of the universal central extensions. Unfortunately this property doesn't
hold for the category of Hom-Leibniz algebras as the following counterexample \ref{counterexample} shows.
This problem lead us to introduce the notion of $\alpha$-central extensions in Definition
\ref{alfacentral}, whose properties relative to the composition are given in Lemma \ref{lema 4}.

\begin{Ex} \label{counterexample}
Consider the two-dimensional  Hom-Leibniz algebra $(L,\alpha_L)$ with basis $\{b_1,b_2\}$, bracket given by $[b_2,b_1]=b_2, [b_2,b_2]=b_1$ (unwritten brackets are equal to zero) and endomorphism $\alpha_L =0$.

Let $(K,\alpha_K)$ be the three-dimensional  Hom-Leibniz algebra with basis $\{a_1, a_2, a_3 \}$, bracket given by $[a_2,a_2]=a_1, [a_3,a_2]=a_3, [a_3,a_3]=a_2$ (unwritten brackets are equal to zero) and  endomorphism $\alpha_K =0$.

Obviously $(K,\alpha_K)$ is a perfect Hom-Leibniz algebra since $K = [K,K]$ and $Z(K) = \langle \{a_1 \} \rangle$.

The linear map $\pi : (K, 0) \to (L,0)$ given by $\pi(a_1)=0, \pi(a_2)=b_1,\pi(a_3)=b_2$, is a  central extension  since  $\pi$ is trivially  surjective and is a homomorphism of Hom-Leibniz algebras:
$$\pi[a_2,a_2]=\pi(a_1)=0; \quad [\pi(a_2),\pi(a_2)]=[b_1,b_1]=0$$
$$\pi[a_3,a_2]=\pi(a_3)=b_2; \quad [\pi(a_3),\pi(a_2)]=[b_2,b_1]=b_2$$
$$\pi[a_3,a_3]=\pi(a_2)=b_1; \quad [\pi(a_3),\pi(a_3)]=[b_2,b_2]=b_1$$

Obviously, $0\cdot \pi = \pi\cdot 0$ and Ker$(\pi) =  \langle \{a_1 \} \rangle$, hence Ker$(\pi) \subseteq Z(K)$.

Now consider  the four-dimensional Hom-Leibniz  algebra $(F, \alpha_F)$ with basis $\{e_1, e_2, e_3, e_4\}$, bracket given by $[e_3, e_3] = e_2, [e_4, e_3] = e_4, [e_4, e_4] = e_3$ (unwritten brackets are equal to zero) and  endomorphism $\alpha_F =0$.

The linear map $\rho : (F, 0) \to (K,0)$ given by $\rho(e_1)=0, \rho(e_2)=a_1,\rho(e_3)=a_2,\rho(e_4)=a_3$, is a  central extension since $\rho$ is trivially   surjective and is a homomorphism of Hom-Leibniz algebras:
$$\rho[e_3,e_2]=\rho(0)= 0; \quad [\rho(e_3),\rho(e_2)]=[a_2,a_1]=0$$
$$\rho[e_3,e_3]=\rho(e_2)= a_1; \quad [\rho(e_3),\rho(e_3)]=[a_2,a_2]=a_1$$
$$\rho[e_4,e_3]=\rho(e_4)= a_3; \quad [\rho(e_4),\rho(e_3)]=[a_3,a_2]=a_3$$
$$\rho[e_4,e_4]=\rho(e_3)= a_2; \quad [\rho(e_4),\rho(e_4)]=[a_3,a_3]=a_2$$
Obviously, $0\cdot \rho = \rho\cdot 0$, Ker$(\rho) =  \langle \{e_1 \} \rangle$ and $Z(F)=  \langle \{e_1 \} \rangle$, hence Ker$(\rho) \subseteq Z(F)$.

The composition $\pi\cdot \rho : (F, 0) \to (L, 0)$ is given by $\pi\cdot \rho(e_1) = \pi(0) = 0,$  $\pi\cdot \rho(e_2) = \pi(a_1) = 0, \pi\cdot \rho(e_3) = \pi(a_2) = b_1, \pi\cdot \rho(e_4) = \pi(a_3) = b_2$. Consequently, $\pi\cdot \rho : (F, 0) \to (L, 0)$ is a surjective homomorphism, but is not a  central extension, since $Z(F) = \langle \{e_1 \} \rangle$ and Ker$(\pi\cdot \rho) = \langle \{ e_1, e_2 \} \rangle$, i.e.  Ker$(\pi\cdot \rho) \nsubseteq Z(F)$.
\end{Ex}

\begin{Le} \label{lema 4}
Let $0 \to (M, \alpha_M) \stackrel{i} \to (K,\alpha_K) \stackrel{\pi} \to (L, \alpha_L) \to 0$  and $0 \to (N, \alpha_N) \stackrel{j} \to (F,\alpha_F) \stackrel{\rho} \to (K, \alpha_K) \to 0$ be central extensions with  $(K, \alpha_K)$ a perfect Hom-Leibniz algebra. Then the composition extension $0 \to (P, \alpha_P) = {\rm Ker}\ (\pi \cdot \rho)   \to (F,\alpha_F) \stackrel{\pi \cdot \rho} \to (L, \alpha_L) \to 0$ is an $\alpha$-central extension.

Moreover, if $0 \to (M, \alpha_M) \stackrel{i} \to (K,\alpha_K) \stackrel{\pi} \to (L, \alpha_L) \to 0$ is a universal $\alpha$-central extension, then $0 \to (N, \alpha_N) \stackrel{j} \to (F,\alpha_F) \stackrel{\rho} \to (K, \alpha_K) \to 0$ is split.
\end{Le}
{\it Proof.} We must proof that $[\alpha_P(P), F] = 0 = [F, \alpha_P(P)]$.

Since $(K,\alpha_K)$ is a perfect Hom-Leibniz algebra, then  every element $f \in F$ can be written as $f = \displaystyle \sum_i \lambda_i [f_{i_1},f_{i_2}]  + n, n \in N, f_{i_j} \in F, j=1,2$. So, for all $p \in P, f \in F$ we have
$$[\alpha_P(p), f] = \displaystyle \sum_i \lambda_i \left( [[p,f_{i_1}],\alpha_F(f_{i_2})] - [[p,f_{i_2}],\alpha_F(f_{i_1})] \right) + [\alpha_P(p),n] = 0$$
since $[p,f_{i_j}] \in {\rm Ker}\ (\rho) \subset Z(F)$.

In a similar way we can check that $[f, \alpha_P(p)]=0$.

For the second statement, if $0 \to (M, \alpha_M) \stackrel{i} \to (K,\alpha_K) \stackrel{\pi} \to (L, \alpha_L) \to 0$ is a universal $\alpha$-central extension, then by the first statement, $0 \to (P, \alpha_P) = {\rm Ker}\ (\pi \cdot \rho)   \to (F,\alpha_F) \stackrel{\pi \cdot \rho} \to (L, \alpha_L) \to 0$ is an $\alpha$-central extension, then there exists a unique homomorphism of Hom-Leibniz algebras $\sigma : (K,\alpha_K) \to (F,\alpha_F)$ such that $\pi \cdot \rho \cdot \sigma = \pi$. On the other hand, $\pi \cdot \rho \cdot \sigma = \pi = \pi \cdot Id$ and $(K,\alpha_K)$ is perfect, then Lemma \ref{lema 2} implies that $\rho \cdot \sigma = Id$. \rdg

\newpage

\begin{Th}\label{teorema}\
\begin{enumerate}
\item[a)]  If a central extension $0 \to (M, \alpha_M)
\stackrel{i} \to (K,\alpha_K) \stackrel{\pi} \to (L, \alpha_L) \to
0$ is a universal $\alpha$-central extension, then
$(K,\alpha_K)$ is a perfect Hom-Leibniz algebra and every central
extension of $(K,\alpha_K)$ splits.

    \item[b)] Let  $0 \to (M, \alpha_M) \stackrel{i} \to (K,\alpha_K) \stackrel{\pi} \to (L, \alpha_L) \to 0$ be a central extension.

     If $(K,\alpha_K)$ is a perfect Hom-Leibniz algebra and every central extension of $(K,\alpha_K)$ splits, then $0 \to (M, \alpha_M) \stackrel{i} \to (K,\alpha_K) \stackrel{\pi} \to (L, \alpha_L) \to 0$ is a universal central extension.

\item[c)] A Hom-Leibniz algebra  $(L, \alpha_L)$ admits a
universal central extension if and only if  $(L, \alpha_L)$ is
perfect.

\item[d)] The kernel of the universal central extension is
canonically isomorphic to $HL_2^{\alpha}(L)$.

\end{enumerate}
\end{Th}
{\it Proof.}\

\noindent {\it a)} If $0 \to (M, \alpha_M) \stackrel{i} \to (K,\alpha_K) \stackrel{\pi} \to (L, \alpha_L) \to 0$ is a universal $\alpha$-central extension, then is a universal central extension by Remark \ref{rem}, so $(K,\alpha_K)$ is a perfect Hom-Leibniz algebra by Lemma \ref{lema 3} and every central extension of $(K,\alpha_K)$ splits by Lemma \ref{lema 4}.

\bigskip

\noindent {\it b)} Let us consider  any central extension $0 \to (N, \alpha_N) \stackrel{j} \to (A,\alpha_A) \stackrel{\tau} \to (L, \alpha_L) \to 0$. Construct the pull-back extension  $0 \to (N, \alpha_N) \stackrel{\chi} \to (P,\alpha_P) \stackrel{\overline{\tau}} \to (K, \alpha_K) \to 0$, where $P = A \times_L K =\{(a,k) \in A \times K \mid \tau(a)=\pi(k) \}$ and $\alpha_P(a,k)=(\alpha_A(a),\alpha_K(k))$, which is central, consequently is split, that is, there exists a homomorphism $\sigma : (K,\alpha_K) \to (P,\alpha_P)$ such that $\overline{\tau} \cdot \sigma = Id$.

Then $\overline{\pi} \cdot \sigma$, where $\overline{\pi} :  (P,\alpha_P)\to (A,\alpha_A)$ is induced by the pull-back construction, satisfies $\tau \cdot\overline{ \pi} \cdot \sigma = \pi$. Lemma \ref{lema 2} ends the proof.

\bigskip

\noindent {\it c)} and {\it d)}   For a Hom-Leibniz algebra $(L, \alpha_L)$  consider the chain homology complex $CL_{\star}^{\alpha}(L)$ which is $CL_{\star}^{\alpha}(L, \mathbb{K})$ where $\mathbb{K}$ is endowed with the trivial Hom-co-representation structure.

As $\mathbb{K}$-vector spaces, let $I_L$ be the subspace of $L \otimes L$ spanned by the elements of the form $-[x_1,x_2] \otimes \alpha_L(x_3) + [x_1,x_3] \otimes \alpha_L(x_2) + \alpha_L(x_1) \otimes [x_2,x_3], x_1, x_2, x_3 \in L$. That is $I_L = {\rm Im}\ \left( d_3 : CL_3^{\alpha}(L) \to CL_2^{\alpha}(L) \right)$.

Now we denote the quotient $\mathbb{K}$-vector space  $\frac{L \otimes L}{I_L}$  by $\frak{uce}(L)$. Every class $x_1 \otimes x_2 + I_L$ is denoted by $\{x_1,x_2\}$, for all $x_1, x_2 \in L$.

By construction, the following identity holds
\begin{equation}
\{\alpha_L(x_1),[x_2,x_3]\} = \{[x_1,x_2],\alpha_L(x_3)\} - \{[x_1,x_3],\alpha_L(x_2)\}
\end{equation}
for all $x_1, x_2, x_3 \in L$.

Now $d_2(I_L)=0$, so it induces a $\mathbb{K}$-linear map $u_L : \frak{uce}(L) \to L$, given by $u_L(\{x_1,x_2\})=[x_1,x_2]$. Moreover $(\frak{uce}(L), \widetilde{\alpha})$, where $\widetilde{\alpha} : \frak{uce}(L) \to \frak{uce}(L)$ is defined by $\widetilde{\alpha}(\{x_1,x_2\}) = \{\alpha_L(x_1), \alpha_L(x_2) \}$,  is a Hom-Leibniz algebra with respect to the bracket $[\{x_1,x_2\},\{y_1,y_2\}]= \{[x_1,x_2],[y_1,y_2]\}$ and $u_L : (\frak{uce}(L), \widetilde{\alpha}) \to (L,\alpha_L)$ is a homomorphism of Hom-Leibniz algebras. Actually, Im $u_L = [L,L]$, but $(L,\alpha_L)$ is a perfect Hom-Leibniz algebra, so $u_L$ is an epimorphism.

From the construction, it follows that Ker $u_L = HL_2^{\alpha}(L)$, so  we have the extension
$$0 \to (HL_2^{\alpha}(L), \widetilde{\alpha}_{\mid}) \to (\frak{uce}(L), \widetilde{\alpha}) \stackrel{u_L}\to (L,\alpha_L) \to 0$$
which is central, because $[{\rm Ker}\ u_L, \frak{uce}(L)] = 0 = [\frak{uce}(L), {\rm Ker}\ u_L]$, and universal because for any central extension $0 \to (M,\alpha_M) \to (K,\alpha_K) \stackrel{\pi} \to (L,\alpha_L) \to 0$ there exists the homomorphism of Hom-Leibniz algebras $\varphi : (\frak{uce}(L), \widetilde{\alpha}) \to (K,\alpha_K)$ given by $\varphi(\{x_1,x_2\})=[k_1,k_2], \pi(k_i)=x_i, i = 1, 2$, such that $\pi \cdot \varphi = u_L$. Moreover, thanks to Lemma \ref{lema 3} and Lemma \ref{lema 2}, $\varphi$ is unique.
\rdg

\begin{Co} \
\begin{enumerate}

 \item[a)] Let $0 \to (M, \alpha_M) \stackrel{i} \to (K,\alpha_K) \stackrel{\pi} \to (L, \alpha_L) \to 0$ be a universal $\alpha$-central extension, then $HL_1^{\alpha}(K) = HL_2^{\alpha}(K) = 0$.

   \item[b)]   Let $0 \to (M, \alpha_M) \stackrel{i} \to (K,\alpha_K) \stackrel{\pi} \to (L, \alpha_L) \to 0$  be a central extension such that $HL_1^{\alpha}(K) = HL_2^{\alpha}(K) = 0$, then $0 \to (M, \alpha_M) \stackrel{i} \to (K,\alpha_K) \stackrel{\pi} \to (L, \alpha_L) \to 0$  is a universal central extension.
     \end{enumerate}
\end{Co}
{\it Proof.}

\noindent {\it a)} If $0 \to (M, \alpha_M) \stackrel{i} \to (K,\alpha_K)
\stackrel{\pi} \to (L, \alpha_L) \to 0$ is a universal
$\alpha$-central extension, then $(K,\alpha_K)$ is perfect by
Remark \ref{rem} and Lemma \ref{lema 3}, so  $HL_1^{\alpha}(K) =
0$. By Lemma \ref{lema 4} and Theorem \ref{teorema} {\it c), d)} the
universal central extension corresponding to $(K,\alpha_K)$ is
split, so $HL_2^{\alpha}(K) = 0$.
\bigskip 

\noindent {\it b)} $HL_1^{\alpha}(K) = 0$ implies that $(K,\alpha_K)$ is a
perfect Hom-Leibniz algebra.

\noindent \quad $HL_2^{\alpha}(K) = 0$ implies that $(\frak{uce}(K),\widetilde{\alpha}) \stackrel{\sim} \to (K,\alpha_K)$. Theorem \ref{teorema} {\it b)} ends the proof. \rdg

\begin{De} \label{central}
An $\alpha$-central extension $0 \to (M, \alpha_M) \stackrel{i} \to (K,\alpha_K) \stackrel{\pi} \to (L, \alpha_L) \to 0$ is said to be universal if for any central extension $0 \to (R, \alpha_R) \stackrel{j} \to (A,\alpha_A) \stackrel{\tau} \to (L, \alpha_L) \to 0$ there exists a unique homomorphism $\varphi : (K,\alpha_K) \to (A, \alpha_A)$ such that $\tau \cdot \varphi = \pi$.

An $\alpha$-central extension $0 \to (M, \alpha_M) \stackrel{i} \to (K,\alpha_K) \stackrel{\pi} \to (L, \alpha_L) \to 0$ is said to be $\alpha$-universal if for any $\alpha$-central extension $0 \to (R, \alpha_R) \stackrel{j} \to (A,\alpha_A) \stackrel{\tau} \to (L, \alpha_L) \to 0$ there exists a unique homomorphism $\psi : (K,\alpha_K) \to (A, \alpha_A)$ such that $\tau \cdot \psi = \pi$.
\end{De}

\begin{Rem}
Obviously, every $\alpha$-universal $\alpha$-central extension is an $\alpha$-central extension which is universal in the sense of Definition \ref{central}. In case $\alpha_M = Id$ both notions coincide with the definition of universal central extension.
\end{Rem}

\begin{Pro}\

\begin{enumerate}
\item[a)] Let $0 \to (M, \alpha_M) \stackrel{i} \to (K,\alpha_K) \stackrel{\pi} \to (L, \alpha_L) \to 0$  and $0 \to (N, \alpha_N) \stackrel{j} \to (F,\alpha_F) \stackrel{\rho} \to (K, \alpha_K) \to 0$ be central extensions. If $0 \to (N, \alpha_N) \stackrel{j} \to (F,\alpha_F) \stackrel{\rho} \to (K, \alpha_K) \to 0$ is a universal central extension, then $0 \to (P, \alpha_P) ={\rm Ker} (\pi \cdot \rho) $ $ \stackrel{\chi}\to (F,\alpha_F) \stackrel{\pi \cdot \rho} \to (L, \alpha_L) \to 0$ is an $\alpha$-central extension which is universal in the sense of Definition \ref{central}.

    \item[b)] Let $0 \to (M, \alpha_M) \stackrel{i} \to (K,\alpha_K) \stackrel{\pi} \to (L, \alpha_L) \to 0$ and $0 \to (N, \alpha_N) \stackrel{j} \to (F,\alpha_F) \stackrel{\rho} \to (K, \alpha_K) \to 0$ be central extensions with $(F,\alpha_F)$ a perfect Hom-Leibniz algebra. If   $0 \to (P, \alpha_P) = {\rm Ker} (\pi \cdot \rho) \stackrel{\chi} \to (F,\alpha_F) \stackrel{\pi \cdot \rho} \to (L, \alpha_L) \to 0$ is an $\alpha$-universal $\alpha$-central extension, then $0 \to (N, \alpha_N) \stackrel{j} \to (F,\alpha_F) \stackrel{\rho} \to (K, \alpha_K) \to 0$ is a universal central extension.
        \end{enumerate}
\end{Pro}
{\it Proof.}
 
 \noindent $a)$ If $0 \to (N, \alpha_N) \stackrel{j} \to (F,\alpha_F) \stackrel{\rho} \to (K, \alpha_K) \to 0$ is a universal central extension, then $(F,\alpha_F)$ and $(K, \alpha_K)$ are perfect Hom-Leibniz algebras by Lemma \ref{lema 3}.

On the other hand, $0 \to (P, \alpha_P) = {\rm Ker} (\pi \cdot \rho) \stackrel{\chi} \to (F,\alpha_F) \stackrel{\pi \cdot \rho} \to (L, \alpha_L) \to 0$ is an $\alpha$-central extension by Lemma \ref{lema 4}.

In order to obtain the universality, for any central extension $0 \to (R, \alpha_R) \to (A,\alpha_A) \stackrel{\tau} \to (L, \alpha_L) \to 0$  construct the pull-back extension corresponding to $\tau$ and $\pi$, $0 \to (R,\alpha_R) \to (K \times_L A, \alpha_K \times \alpha_A) \stackrel{\overline{\tau}}\to (K, \alpha_K) \to 0$. Since $0 \to (N, \alpha_N) \stackrel{j} \to (F,\alpha_F) \stackrel{\rho} \to (K, \alpha_K) \to 0$ is a universal central extension, then there exists a unique homomorphism $\varphi : (F, \alpha_F) \to (K \times_L A, \alpha_K \times \alpha_A)$ such that $\overline{\tau} \cdot \varphi = \rho$. Then the homomorphism $\overline{\pi} \cdot \varphi$ satisfies that $\tau \cdot \overline{\pi} \cdot \varphi = \pi \cdot \rho$ and it is unique by Lemma \ref{lema 2}.
\bigskip

\noindent $b)$ $(F, \alpha_F)$ perfect implies that $(K,\alpha_K)$ is perfect by Lemma \ref{lema 1}. In order to prove the universality of the central extension $0 \to (N, \alpha_N) \stackrel{j} \to (F,\alpha_F) \stackrel{\rho} \to (K, \alpha_K) \to 0$, let us consider any central extension $0 \to (R, \alpha_R) \to (A,\alpha_A) \stackrel{\sigma} \to (K, \alpha_K ) \to 0$, then $0 \to {\rm Ker}(\pi \cdot \sigma) \to (A, \alpha_A)  \stackrel{\pi \cdot \sigma}\to (L, \alpha_L) \to 0$ is an $\alpha$-central extension by Lemma \ref{lema 4}.

The $\alpha$-universality of $0 \to (P, \alpha_P) = {\rm Ker} (\pi \cdot \rho) \stackrel{\chi} \to (F,\alpha_F) \stackrel{\pi \cdot \rho} \to (L, \alpha_L) \to 0$ implies the existence of a unique homomorphism $\omega : (F, \alpha_F) \to (A, \alpha_A)$ such that $\pi \cdot \sigma \cdot \omega = \pi \cdot \rho$.

Since $(F,\alpha_F)$ is a perfect Hom-Leibniz algebra and Im$(\sigma \cdot \omega) \subseteq$ Im$(\rho)$ + Ker$(\pi)$, then easily is followed that $\sigma \cdot \omega = \rho$. The uniqueness is followed from Lemma \ref{lema 2}. \rdg

\section{Relationships between   universal $\alpha$-central extensions}

Since the  universal central extensions of a perfect  Lie algebra in the categories of Lie and Leibniz algebras are related by means of the results given in  \cite{Gn 1}, then it is natural to consider the generalization of these results to the framework of Hom-Leibniz algebras. But, once more, through the proofs it is necessary to deal with the composition of central extensions, that is, we need to apply   Lemma \ref{lema 4}. This fact causes no significant relationships between the universal central extensions of a perfect Hom-Lie algebra in the categories {\sf Hom-Lie} and {\sf Hom-Leib}, so it is not possible to obtain a complete generalization of the mentioned results in \cite{Gn 1}, but it is possible to obtain results in the more restricted context of universal $\alpha$-central extensions.

We start with the introduction of a new concept of perfection.

 \begin{De}
 A Hom-Lie algebra (respectively, Hom-Leibniz) $\left(  L, \alpha_{L}\right)$ is said to be $\alpha$-perfect if
 $$L = [\alpha_L(L), \alpha_L(L)]$$
 \end{De}

\begin{Ex}
 The three-dimensional Hom-Lie algebra  $(L,[-,-],\alpha_L)$ with basis $\{a_1,a_2,a_3\}$, bracket given by $[a_1,a_2]= - [a_2,a_1] = a_3; [a_2,a_3]= - [a_3,a_2] = a_1; [a_3,a_1]= - [a_1,a_3] = a_2$ (unwritten brackets are equal to zero) and endomorphism $\alpha_L$ represented by the matrix $\left( \begin{array}{ccc} \frac{\sqrt{2}}{2} & 0 & \frac{\sqrt{2}}{2} \\ 0 & -1 & 0 \\ \frac{\sqrt{2}}{2} & 0 & - \frac{\sqrt{2}}{2} \end{array} \right)$ is an $\alpha$-perfect Hom-Lie algebra.
\end{Ex}

 \begin{Rem}\label{alfa perfecta} \
 \begin{enumerate}
 \item[a)] When $\alpha_L= Id$, the notion of $\alpha$-perfection coincides with the  notion of perfection.

 \item[b)] Obviously, if $\left(  L, [-,-],\alpha_{L}\right)$ is an $\alpha$-perfect Hom-Lie algebra (respectively, Hom-Leibniz), then it is perfect. Nevertheless the converse is not true in general. For example, the three-dimensional Hom-Lie  algebra $\left(  L, [-,-],\alpha_{L}\right)$  with basis $\{a_1, a_2, a_3\}$, bracket given by $[a_1,a_2]=-[a_2,a_1]=a_3; [a_1,a_3]=-[a_3,a_1]=a_2; [a_2,a_3]=-[a_3,a_2]=a_1$ and endomorphism $\alpha_L=0$ is perfect, but it is not $\alpha$-perfect.

     \item[c)] If  $\left(  L, [-,-],\alpha_{L}\right)$  is $\alpha$-perfect, then $L = \alpha_L(L)$, i.e. $\alpha_L$ is surjective. Nevertheless the converse is not true. For instance, the two-dimensional Hom-Lie algebra with basis $\{a_1, a_2\}$, bracket given by $[a_1, a_2] = - [a_2, a_1] = a_2$ and endomorphism $\alpha_L$ represented by the matrix $\left( \begin{array}{cc} 1 & 0 \\ 0 & 2 \end{array} \right)$. Obviously the endomorphism $\alpha_L$ is surjective, but $[\alpha_L(L), \alpha_L(L)] = \langle \{a_2 \} \rangle$.
 \end{enumerate}
 \end{Rem}

 \begin{Le} \label{lema 10}
Let $0 \to (M, \alpha_M) \stackrel{i} \to (K,\alpha_K) \stackrel{\pi} \to (L, \alpha_L) \to 0$ be a central extension and $(K,\alpha_K)$ an $\alpha$-perfect Hom-Lie (Hom-Leibniz) algebra. If there exists a homomorphism of Hom-Lie (Hom-Leibniz) algebras $f : (K,\alpha_K) \to (A, \alpha_A)$ such that $\tau \cdot f = \pi$, where $0 \to (N, \alpha_N) \stackrel{j} \to (A,\alpha_A) \stackrel{\tau} \to (L, \alpha_L) \to 0$ is an $\alpha$-central extension, then $f$ is unique.
\end{Le}
{\it Proof.} Let us assume that there exist two homomorphisms  $f_1$ and $f_2$ such that $\tau \cdot f_1 = \pi = \tau \cdot f_2$, then $f_1 - f_2 \in$ Ker $\tau = N$, i. e. $f_1(k) = f_2(k)+n_k, n_k \in N$, and $\alpha_N(N) \subseteq Z(A)$.

Moreover, $f_1$ and $f_2$ coincide over $[\alpha_K(K),\alpha_K(K)]$. Indeed $f_1[\alpha_K(k_1),\alpha_K(k_2)]$ $=[\alpha_A \cdot f_1(k_1), \alpha_A \cdot f_1(k_2)] =[\alpha_A \cdot f_2(k_1) + \alpha_A(n_{k_1}), \alpha_A \cdot f_2(k_2) + \alpha_A(n_{k_2})] = [\alpha_A \cdot f_2(k_1) , \alpha_A \cdot f_2(k_2)] = f_2[\alpha_K(k_1),\alpha_K(k_2)]$.

Since $(K,\alpha_K)$ is $\alpha$-perfect, then $f_1$ coincides with $f_2$ over $K$. \rdg

\begin{Th}\label{alfa uce} \
\begin{enumerate}
\item[a)] An $\alpha$-perfect Hom-Lie algebra admits a universal $\alpha$-central extension.
\item[b)] An $\alpha$-perfect Hom-Leibniz algebra admits a universal $\alpha$-central extension.
\end{enumerate}
\end{Th}
{\it Proof.} 

\noindent $a)$ For the Hom-Lie algebra $(L,[-,-],\alpha_L)$, consider the chain complex $$C_{\star}^{\alpha}(L, \mathbb{K}) :  ( C_n^{\alpha}(L,\mathbb{K}) = \Lambda^n L, d_n: C_n^{\alpha}(L,\mathbb{K}) \to C_{n-1}^{\alpha}(L,\mathbb{K}), $$ $$d_n(x_1 \wedge \dots \wedge x_n ) =\displaystyle \sum_{1 \leq i < j \leq n} (-1)^{i+j} [x_i,x_j] \wedge \alpha_L(x_1) \wedge \dots \wedge \widehat{\alpha_L(x_i)} \wedge \dots \wedge \widehat{\alpha_(x_j)} \wedge \dots \alpha_L(x_n)),$$ where $\mathbb{K}$ is considered with a trivial Hom-L-module structure.

Now we consider the quotient $\mathbb{K}$-vector space  $\frak{uce}^{\rm Lie}_{\alpha}(L) = \frac{\alpha_L(L) \wedge \alpha_L(L)}{I_L}$, where $I_L$ is the vector subspace of $\alpha_L(L) \wedge \alpha_L(L)$ spanned by the elements of the form
$$-[x_1,x_2] \wedge \alpha_L(x_3) + [x_1,x_3] \wedge \alpha_L(x_2) -  [x_2,x_3] \wedge \alpha_L(x_1)$$
for all $x_1, x_2, x_3 \in L$.

Observe that every summand of the form $[x_1,x_2] \wedge \alpha_L(x_3)$  is an element of $\alpha_L(L) \wedge \alpha_L(L)$, since  $L$  is $\alpha$-perfect and then $[x_1,x_2] \in L = [\alpha_L(L), \alpha_L(L)] \subseteq \alpha_L(L)$.

We denote by $\{\alpha_L(x_1), \alpha_L(x_2)\}$ the equivalence class $\alpha_L(x_1) \wedge \alpha_L(x_2) + I_L$. $\frak{uce}^{\rm Lie}_{\alpha}(L)$ is endowed with a structure  of Hom-Lie algebra with respect to the bracket $$[\{\alpha_L(x_1),\alpha_L(x_2)\}, \{\alpha_L(y_1),\alpha_L(y_2)\}] = \{ [\alpha_L(x_1),\alpha_L(x_2)], [\alpha_L(y_1),\alpha_L(y_2)]\}$$
and the endomorphism $\widetilde{\alpha} : \frak{uce}^{\rm Lie}_{\alpha}(L) \to \frak{uce}^{\rm Lie}_{\alpha}(L)$ given by $\widetilde{\alpha}(\{\alpha_L(x_1),\alpha_L(x_2)\}) = \{\alpha_L^2(x_1),\alpha_L^2(x_2)\}$.

The restriction of $d_2 : C_2^{\alpha}(L) = L \wedge L \to C_1^{\alpha}(L) = L$ to $\alpha_L(L) \wedge \alpha_L(L)$ vanishes on $I_L$, then it  induces a linear map $u_{\alpha} : \frak{uce}^{\rm Lie}_{\alpha}(L) \to L$, which is defined by $u_{\alpha}(\{\alpha_L(x_1), \alpha_L(x_2)\})= [\alpha_L(x_1), \alpha_L(x_2)]$. Moreover, since $(L,[-,-],\alpha_L)$ is $\alpha$-perfect, then $u_{\alpha}$ is a surjective homomorphism, because ${\rm Im}(u_{\alpha}) =[\alpha_L(L) ,\alpha_L(L)] = L$ and
\begin{enumerate}
\item[$\centerdot$]
$u_{\alpha} [\{\alpha_L(x_1), \alpha_L(x_2) \},\{\alpha_L(y_1), \alpha_L(y_2) \}] = u_{\alpha} \{ [\alpha_L(x_1), \alpha_L(x_2)],[\alpha_L(y_1), \alpha_L(y_2)] \}  =[ [\alpha_L(x_1), \alpha_L(x_2)],[\alpha_L(y_1), \alpha_L(y_2)]]$ $ = [u_{\alpha}\{\alpha_L(x_1), \alpha_L(x_2) \},u_{\alpha} \{\alpha_L(y_1), \alpha_L(y_2) \}].$

\item[$\centerdot \centerdot$] $\alpha_L \cdot u_{\alpha} \{\alpha_L(x_1), \alpha_L(x_2) \} \  = \ \alpha_L [\alpha_L(x_1), \alpha_L(x_2)] \ = \ [\alpha_L^2(x_1), \alpha_L^2(x_2)]=$
$u_{\alpha} \{\alpha_L^2(x_1), \alpha_L^2(x_2)\}\ = \ u_{\alpha} \cdot \widetilde{\alpha}(\{\alpha_L(x_1), \alpha_L(x_2)\}).$
\end{enumerate}
Then we have constructed the  extension
$$0 \to (Ker (u_{\alpha}), \widetilde{\alpha}_{\mid}) \to (\frak{uce}_{\alpha}^{\rm Lie}(L), \widetilde{\alpha})  \stackrel{u_{\alpha}} \to (L, \alpha_L) \to 0$$
which is a central extension. Indeed, for any  $\{\alpha_L(x_1),
\alpha_L(x_2)\} \in Ker (u_{\alpha}),$ $\{\alpha_L(y_1),
\alpha_L(y_2)\} \in \frak{uce}^{\rm Lie}_{\alpha}(L)$, one
verifies
$$[\{\alpha_L(x_1), \alpha_L(x_2)\},\{\alpha_L(y_1),
\alpha_L(y_2)\}] = \{[\alpha_L(x_1),
\alpha_L(x_2)],[\alpha_L(y_1), \alpha_L(y_2)]\} =0$$ since
$\{\alpha_L(x_1), \alpha_L(x_2)\} \in  Ker (u_{\alpha})
\Longleftrightarrow [\alpha_L(x_1), \alpha_L(x_2)] = 0$.

Finally, this central extension is universal $\alpha$-central.
Indeed, consider an $\alpha$-central extension $0 \to (M,
\alpha_M) \stackrel{j}\to (K, \alpha_K) \stackrel{\pi} \to (L,
\alpha_L) \to 0$. We define the map $\Phi : (\frak{uce}^{\rm
Lie}_{\alpha}(L), \widetilde{\alpha}) \to (K, \alpha_K)$ by
$\Phi(\{\alpha_L(x_1), \alpha_L(x_2)\})= [\alpha_K(k_1),
\alpha_K(k_2)]$, where $\pi(k_i)=x_i, x_i \in L, k_i \in K,
i=1,2$. Now we are going to check that $\Phi$ is a homomorphism of
 Hom-Lie algebras such that $\pi \cdot \Phi = u_{\alpha}$. Indeed:
\begin{enumerate}
\item[$\centerdot$] Well-definition: if  $\pi(k_i)=x_i=\pi(k_i'),
i = 1, 2$, then $k_i - k_i' \in Ker (\pi) = M$, so $\alpha_K(k_i)
- \alpha_K(k_i') \in \alpha_M(M)$. Hence,
    $\Phi(\{\alpha_L(x_1), \alpha_L(x_2)\}) = [\alpha_K(k_1), \alpha_K(k_2)]$ $= [\alpha_K(k_1')+ \alpha_M(m_1), \alpha_L(k_2')+ \alpha_M(m_2)] = [\alpha_K(k_1'),\alpha_L(k_2')]$.

    Observe that the condition of   $\alpha$-centrality is essential at this step  in order to be guaranteed
     that $[\alpha_M(M) , K ]= 0$.
\item[$\centerdot \centerdot$]  $\Phi[\{\alpha_L(x_1),
\alpha_L(x_2)\},\{\alpha_L(y_1), \alpha_L(y_2)\}] =
\Phi\{[\alpha_L(x_1), \alpha_L(x_2)],[\alpha_L(y_1),
\alpha_L(y_2)]\}$ $= [[\alpha_L(k_1), \alpha_L(k_2)],[\alpha_L(l_1),
\alpha_L(l_2)]] = [\Phi\{\alpha_L(x_1),
\alpha_L(x_2)\},\Phi\{\alpha_L(y_1), \alpha_L(y_2)\}]$

 where $\pi(k_i)=x_i, \pi(y_i)=l_i, i=1,2$.

 \item[$\centerdot \centerdot \centerdot$] $\alpha_K \cdot \Phi(\{\alpha_L(x_1), \alpha_L(x_2)\}) =
 \alpha_K[\{\alpha_K(k_1), \alpha_K(k_2)\}] = [\alpha_K^2(k_1),
 \alpha_K^2(k_2)] = \Phi \{\alpha_L^2(x_1), \alpha_L^2(x_2)\}  =
 \Phi \cdot \widetilde{\alpha} \{\alpha_L(x_1), \alpha_L(x_2)\}$.

     \item[$\centerdot \centerdot \centerdot \centerdot$] $\pi \cdot \Phi \{\alpha_L(x_1), \alpha_L(x_2)\} = \pi [\alpha_K(k_1), \alpha_K(k_2)] = [\pi \cdot \alpha_K(k_1), \pi \cdot \alpha_K(k_2)] = [\alpha_L \cdot \pi(k_1), \alpha_L \cdot \pi(k_2)] = [\alpha_L(x_1), \alpha_L(x_2)] = u_{\alpha} \{\alpha_L(x_1), \alpha_L(x_2)\}$.
\end{enumerate}
To end the proof, we must verify the uniqueness of $\Phi$: first
at all, we check that $\frak{uce}^{\rm Lie}_{\alpha}(L)$ is
$\alpha$-perfect. Indeed,
 $$[\widetilde{\alpha} \{\alpha_L(x_1), \alpha_L(x_2)\}, \widetilde{\alpha} \{\alpha_L(y_1), \alpha_L(y_2)\}] = \{[\alpha_L^2(x_1), \alpha_L^2(x_2)],[\alpha_L^2(y_1), \alpha_L^2(y_2)]\}$$
i.e. $[\widetilde{\alpha} (\frak{uce}^{\rm Lie}_{\alpha}(L)),
\widetilde{\alpha}( \frak{uce}^{\rm Lie}_{\alpha}(L))] \subseteq
\frak{uce}^{\rm Lie}_{\alpha}(L)$.

For the converse inclusion, having in mind that $\alpha_L(x_i) \in
L = [\alpha_L(L), \alpha_L(L)]$, we have that:
 $$\{\alpha_L(x_1), \alpha_L(x_2)\} = \left \{\alpha \left( \sum_i \lambda_i [\alpha_L(l_{i_1}), \alpha_L(l_{i_2})] \right), \alpha \left( \sum_j \mu_j [\alpha_L(l'_{j_1}), \alpha_L(l'_{j_2})] \right)\right \}$$
 $$\sum \lambda_i \mu_j \left\{ [\alpha^2_L(l_{i_1}),\alpha^2_L(l_{i_2})], [\alpha^2_L(l'_{j_1}),\alpha^2_L(l'_{j_2})] \right \} =$$
 $$\sum \lambda_i \mu_j \left[ \{\alpha^2_L(l_{i_1}),\alpha^2_L(l_{i_2})\} , \{\alpha^2_L(l'_{j_1}),\alpha^2_L(l'_{j_2})\} \right] =$$
$$\sum \lambda_i \mu_j \left[ \widetilde{\alpha} \{\alpha_L(l_{i_1}),\alpha_L(l_{i_2})\}, \widetilde{\alpha} \{\alpha_L(l'_{j_1}),\alpha_L(l'_{j_2})\} \right] \in [\widetilde{\alpha} (\frak{uce}^{\rm Lie}_{\alpha}(L)),
\widetilde{\alpha}( \frak{uce}^{\rm Lie}_{\alpha}(L))]$$

Now Lemma \ref{lema 10} ends the proof.
\bigskip

\noindent $b)$ For the Hom-Leibniz algebra $(L,[-,-],\alpha_L)$, we consider
the chain complex $CL_{\star}^{\alpha}(L, \mathbb{K})$, where
$\mathbb{K}$ is considered with trivial  Hom-co-representation
structure.

Take the quotient vector space $\frak{uce}^{\rm Leib}_{\alpha}(L)
= \frac{\alpha_L(L) \otimes \alpha_L(L)}{I_L}$, where $I_L$ is the
vector  subspace of $\alpha_L(L) \otimes \alpha_L(L)$ generated by the
elements of the form
$$-[x_1,x_2] \otimes \alpha_L(x_3) + [x_1,x_3] \otimes \alpha_L(x_2) + \alpha_L(x_1) \otimes [x_2,x_3]$$
for all $x_1, x_2, x_3 \in L$.

Observe that every summand of the form $[x_1,x_2] \otimes
\alpha_L(x_3)$ or $\alpha_L(x_1) \otimes [x_2,x_3]$ belongs to
$\alpha_L(L) \otimes \alpha_L(L)$, since $L$ is $\alpha$-perfect,
then $[x_1,x_2] \in L = [\alpha_L(L), \alpha_L(L)] \subseteq
\alpha_L(L)$.

We denote by $\{\alpha_L(x_1),
\alpha_L(x_2)\}$ the equivalence class $\alpha_L(x_1) \otimes
\alpha_L(x_2) + I_L$. $\frak{uce}^{\rm Leib}_{\alpha}(L)$ is
endowed with a structure of Hom-Leibniz algebra with respect to
the bracket
$$[\{\alpha_L(x_1),\alpha_L(x_2)\},
\{\alpha_L(y_1),\alpha_L(y_2)\}] = \{
[\alpha_L(x_1),\alpha_L(x_2)], [\alpha_L(y_1),\alpha_L(y_2)]\}$$
and the endomorphism $\overline{\alpha} : \frak{uce}^{\rm
Leib}_{\alpha}(L) \to \frak{uce}^{\rm Leib}_{\alpha}(L)$ defined
by $\overline{\alpha}(\{\alpha_L(x_1),\alpha_L(x_2)\})$ $=
\{\alpha_L^2(x_1),\alpha_L^2(x_2)\}$.

The restriction of the differential $d_2 : CL_2^{\alpha}(L) = L
\otimes L \to CL_1^{\alpha}(L) = L$ to $\alpha_L(L) \otimes
\alpha_L(L)$ vanishes on $I_L$, so it induces a linear map
$U_{\alpha} : \frak{uce}^{\rm Leib}_{\alpha}(L) \to L$, that is
given by $U_{\alpha}(\{\alpha_L(x_1), \alpha_L(x_2)\})=
[\alpha_L(x_1), \alpha_L(x_2)]$. Moreover, thanks to be
$(L,[-,-],\alpha_L)$  $\alpha$-perfect, then $U_{\alpha}$ is a
surjective homomorphism, since Im$(U_{\alpha}) =[\alpha_L(L)
,\alpha_L(L)]=L$ and
\begin{enumerate}
\item[$\centerdot$] $U_{\alpha} [\{\alpha_L(x_1), \alpha_L(x_2)
\},\{\alpha_L(y_1), \alpha_L(y_2) \}] = U_{\alpha} \{
[\alpha_L(x_1), \alpha_L(x_2)],[\alpha_L(y_1), \alpha_L(y_2)] \} =
[ [\alpha_L(x_1), \alpha_L(x_2)],[\alpha_L(y_1), \alpha_L(y_2)] ]=
[U_{\alpha}\{\alpha_L(x_1), \alpha_L(x_2) \},U_{\alpha}
\{\alpha_L(y_1), \alpha_L(y_2) \}].$

\item[$\centerdot \centerdot$] $\alpha_L \cdot U_{\alpha}
\{\alpha_L(x_1), \alpha_L(x_2) \}\ =\ \alpha_L [\alpha_L(x_1),
\alpha_L(x_2)] \ = \ [\alpha_L^2(x_1), \alpha_L^2(x_2)] = U_{\alpha}
\{\alpha_L^2(x_1), \alpha_L^2(x_2)\} = U_{\alpha} \cdot
\overline{\alpha} \{\alpha_L(x_1), \alpha_L(x_2) \}.$
\end{enumerate}
so  we have  constructed the extension
$$0 \to (Ker (U_{\alpha}), \overline{\alpha}_{\mid}) \to (\frak{uce}_{\alpha}^{\rm Leib}(L), \overline{\alpha})  \stackrel{U_{\alpha}} \to (L, \alpha_L) \to 0$$
which is central. Indeed, for any $\{\alpha_L(x_1),
\alpha_L(x_2)\} \in Ker (U_{\alpha}), \{\alpha_L(y_1),
\alpha_L(y_2)\} \in \frak{uce}_{\alpha}^{\rm Leib}(L)$, we have that
$$[\{\alpha_L(x_1), \alpha_L(x_2)\},\{\alpha_L(y_1),
\alpha_L(y_2)\}] = \{[\alpha_L(x_1),
\alpha_L(x_2)],[\alpha_L(y_1), \alpha_L(y_2)]\} =0$$ since
$\{\alpha_L(x_1), \alpha_L(x_2)\} \in  Ker (U_{\alpha})
\Longleftrightarrow [\alpha_L(x_1), \alpha_L(x_2)] = 0$.

In a similar way it is verified that $[\frak{uce}^{\rm Leib}_{\alpha}(L), Ker
(U_{\alpha})] = 0$.

The checking that this central extension is universal
$\alpha$-central  follows parallel arguments as in Hom-Lie case,
so we omit it. \rdg

\bigskip

Let $\left(  L,\alpha_{L}\right)$ be an $\alpha$-perfect Hom-Lie
algebra. By Theorem \ref{alfa uce}, $\left( L,\alpha_{L}\right) $
admites a universal $\alpha$-central extension $\left(
\frak{uce}_{\alpha}^{\rm Lie}\left(  L\right),
\widetilde{\alpha} \right)$ in the category {\sf Hom-Lie} of Hom-Lie algebras  and a universal
$\alpha$-central extension $\left(
\frak{uce}_{\alpha}^{\rm Leib}\left( L\right),\overline{ \alpha} \right)$ in the category {\sf Hom-Leibniz}. The
following result provides a relationship between both extensions.

\begin{Pro}
$\left(  \frak{uce}_{\alpha}^{\rm Leib}\left( L\right), \overline{\alpha} \right)$ is the universal central extension of
the $\alpha$-perfect Hom-Lie algebra $\left( \frak{uce}_{\alpha}^{\rm Lie}\left( L\right), \widetilde{\alpha} \right)$ in the category {\sf Hom-Leib}.

\noindent Moreover there is an isomorphism of Hom-Lie algebras
$$\left(  \frak{uce}_{\alpha}^{\rm Lie}\left(  L\right) , \widetilde{\alpha}
\right)  \cong \left(  \left( \frak{uce}_{\alpha}^{\rm Leib}\left(  L\right) \right)_{\rm Lie}
,\overline{\alpha}_{\rm Lie}\right).$$
\end{Pro}
{\it Proof.} Since $\left(  L,\alpha_{L}\right)$ is an
$\alpha$-perfect Hom-Lie algebra, then by Theorem \ref{alfa uce} it admits a universal $\alpha$-central extension in the categories {\sf Hom-Lie} and {\sf Hom-Leib}. This situation is described in the following diagram:
$$\xymatrix{
  0  \ar[r] & \left(  Ker\left(  U_{\alpha}\right), \overline{\alpha_{\mid}}\right) \ar[r]  & \left(  \frak{uce}_{\alpha}^{\rm Leib}\left(  L\right),\overline{\alpha} \right) \ar[r]^{\ \ \ U_{\alpha}} \ar@{-->}[d]^{\exists! \Phi}& \left(  L,\alpha_{L}\right) \ar[r]\ar@{=}[d] & 0 \quad ({\sf Hom-Leib})\\
 0 \ar[r] &  \left(  Ker\left(  u_{\alpha}\right), \widetilde{\alpha_{\mid}}\right) \ar[r] & \left(  \frak{uce}_{\alpha}^{\rm Lie}\left(  L\right), \widetilde{\alpha}\right) \ar[r]^{\ \ \ u_{\alpha}} & \left(  L,\alpha_{L}\right) \ar[r] & 0 \quad  ({\sf Hom-Lie})}$$
Since $U_{\alpha} : \left(  \frak{uce}_{\alpha}^{\rm Leib}\left(
L\right) ,\overline{\alpha} \right) \to \left(
L,\alpha_{L}\right)$ is a universal $\alpha$-central extension, by
Remark \ref{rem}, it is a universal central extension; since
$u_{\alpha} : \left( \frak{uce}_{\alpha}^{\rm Lie}\left(  L\right)
,\widetilde{\alpha}\right) \to \left(  L,\alpha_{L}\right)$ is a
 central extension, then there exists a unique homomorphism of Hom-Leibniz algebras $\Phi : \left(
\frak{uce}_{\alpha}^{\rm Leib}\left(  L\right)
,\overline{\alpha} \right) \to \left(
\frak{uce}_{\alpha}^{\rm Lie}\left(  L\right)
,\widetilde{\alpha} \right)$ such that
 $u_{\alpha}.\Phi=U_{\alpha}$.

We only need to prove that the extension
$$0 \to  \left(  Ker\left(  \Phi\right), \overline{\alpha_{\mid}} \right)  \to \left(  \frak{uce}_{\alpha}^{\rm Leib}\left(  L\right),\overline{\alpha} \right) \stackrel{\Phi}\to \left(
\frak{uce}_{\alpha}^{\rm Lie}\left(  L\right),\widetilde{\alpha} \right)
\to 0$$ is a universal central extension in the category {\sf
Hom-Leib}.

First at all, by construction, we have that $Ker(\Phi) \subseteq
Ker(U_{\alpha}) \subseteq Z\left( \frak{uce}_{\alpha}^{\rm Leib}\left(
L\right)\right)$.

On the other hand,  $\Phi$ is a surjective homomorphism, since
$\frak{uce}_{\alpha}^{\rm Lie}\left(  L\right) \subseteq
\operatorname{Im}\left(  \Phi\right) +Ker\left(
u_{\alpha}\right)$. Indeed, for $z' \in \frak{uce}_{\alpha}^{\rm Lie}
\left(  L\right), u_{\alpha}(z') \in L$, then there exists any
$z\in \frak{uce}_{\alpha}^{\rm Leib}\left( L\right)$ such that
$U_{\alpha}\left(  z\right) =u_{\alpha}\left(z'\right)$, but
$U_{\alpha} =u_{\alpha}\cdot  \Phi$, then $z' - \Phi\left(z\right) \in
Ker\left( u_{\alpha}\right)$ and,  consequently, $z' \in Im \left(
\Phi\right)  + Ker\left( u_{\alpha}\right)$.

Since $u_{\alpha} : \left(  \frak{uce}_{\alpha}^{\rm Lie}\left(
L\right),\widetilde{\alpha} \right) \to \left(  L,\alpha_{L}\right)$
is a universal  $\alpha$-central extension, then by Remark \ref{rem} in \cite{CIP} is a
universal central extensions and, by Lemma \ref{lema 3} in \cite{CIP}, we know
that $\left( \frak{uce}_{\alpha}^{\rm Lie}\left(
L\right),\widetilde{\alpha} \right)$ is perfect. Hence:
$$\frak{uce}_{\alpha}^{\rm Lie}\left(  L\right)  = \left[  \frak{uce}_{\alpha
}^{\rm Lie}\left(  L\right), \frak{uce}_{\alpha}^{\rm Lie}\left(  L\right)  \right]
\subseteq$$
$$\left[  Im \left(  \Phi\right)  +Ker\left(  u_{\alpha}
\right),Im\left(  \Phi\right)  +Ker\left(  u_{\alpha }\right)
\right]  \subseteq Im \left(  \Phi\right)$$ since $\left[ Im\left(
\Phi\right), Ker\left(  u_{\alpha }\right)  \right]  = 0 = \left[
Ker\left(  u_{\alpha}\right), Im\left(  \Phi\right) \right]$.

Now we check that the conditions of Theorem \ref{teorema}  {\it b)} hold:
\begin{enumerate}
\item[$\centerdot$] Since $U_{\alpha} :
\left(\frak{uce}_{\alpha}^{\rm Leib}\left( L\right),
\overline{\alpha} \right) \to \left(  L,\alpha_{L}\right)$
is a universal $\alpha$-central extension, then by Remark \ref{rem} is a universal central extension and, by Lemma \ref{lema 3},
 we know that  $\left(  \frak{uce}_{\alpha}^{\rm Leib}\left(  L\right),\overline{\alpha} \right)$ is perfect.

\item[$\centerdot \centerdot$] Now we must check that any cental
extension of the form $0 \to \left( P,\alpha_{P}\right) \to \left(
K,\alpha_{K}\right)\overset{\rho}{\to}\left(
\frak{uce}_{\alpha}^{\rm Leib}\left(
L\right),\overline{\alpha} \right)$ is split.

 For that, we consider the composition extension $U_{\alpha}\cdot \rho : \left(  K,\alpha_{K}\right)  \to \left(
L,\alpha_{L}\right)$, which is an $\alpha-$central extension
thanks to Lemma \ref{lema 4}. Now Definition \ref{Def universal}
guaranties the existence of a unique homomorphism of Hom-Leibniz
algebras $\sigma:\left( \frak{uce}_{\alpha}^{\rm Leib}\left(
L\right),\overline{\alpha}\right)\to \left(  \
K,\alpha_{K}\right)$ such that $U_{\alpha}\cdot \rho \cdot \sigma =
U_{\alpha}$.

Since $U_{\alpha} \cdot Id = U_{\alpha}$, then Lemma \ref{lema 2}
implies that $\rho \cdot \sigma = Id$.
\end{enumerate}

In order to prove the isomorphism established in the second
statement, we  consider the following diagram where $\overline{\Phi}$ is
induced by $\Phi$ and $\gamma$ by $U_{\alpha}$:
\[\xymatrix{& & \left(  \frak{uce}_{\alpha}^{\rm Leib}\left(  L\right),\overline{\alpha} \right)^{\rm ann} \ar@{>->}[dl]&
\\
 & \left(  \frak{uce}_{\alpha}^{\rm Leib}\left(  L\right),\overline{\alpha} \right) \ar[r]^{\ \ U_{\alpha}} \ar@{-}[d] \ar@{>>}[dl]& \left(  L,\alpha_{L}\right) \ar[r]\ar@{=}[d] & 0 \\
  \left(  \left(\frak{uce}_{\alpha}^{\rm Leib}\left(  L\right) \right)_{\rm Lie} ,\overline{\alpha}_{\rm Lie}\right) \ar@{-->>}^{\overline{\Phi}}[rd] \ar@{-->>}^{\gamma}[rru]& \ar@{>>}[d]^{\Phi} & \ar@{=}[d]\\
&  \left(  \frak{uce}_{\alpha}^{\rm Lie}\left(  L\right)
,\widetilde{\alpha} \right) \ar[r]^{\ \ u_{\alpha}} & \left(
L,\alpha_{L}\right) \ar[r] & 0   } \] 

$\gamma : \left(  \left(\frak{uce}_{\alpha}^{\rm Leib}\left(  L\right)
\right)_{\rm Lie} ,\overline{\alpha}_{\rm Lie}\right) \to \left(
L,\alpha_{L}\right)$ is a  central extension since $\gamma$ is a
surjective homomorphism. $Ker\left( \gamma\right) \subseteq
Z\left( \left( \frak{uce}_{\alpha}^{\rm Leib}\left( L\right) \right)
_{\rm Lie}\right)$ since for $\overline{x}\in Ker\left( \gamma
\right)$ and $\overline{y}\in\left(
\frak{uce}_{\alpha}^{\rm Leib}\left( L\right) \right)  _{\rm Lie}$,
$\left[ \overline{x},\overline{y}\right] =\overline{\left[
x,y\right] }=\overline{0}$, because $0 = \gamma(\overline{x}) =
U_{\alpha}(x)$, which implies that $x \in Ker(U_{\alpha})
\subseteq Z(\frak{uce}_{\alpha}^{\rm Leib}\left( L\right) )$.

As $u_{\alpha} : \left(  \frak{uce}_{\alpha}^{\rm Lie}\left(  L\right)
,\widetilde{\alpha} \right)\to \left(  L,\alpha_{L}\right)$ is a
universal $\alpha$-central extension, then is a universal central
extension by Remark \ref{rem} in \cite{CIP}, and since $\gamma : \left(
\left(\frak{uce}_{\alpha}^{\rm Leib}\left( L\right) \right)_{\rm Lie}
,\overline{\alpha}_{\rm Lie}\right)$ $\to \left(
L,\alpha_{L}\right)$ is a central extension, then there exists a
unique homomorphism $\Psi : \left( \frak{uce}_{\alpha}^{\rm Lie}\left(
L\right) ,\widetilde{\alpha} \right) \to \left(
\left(\frak{uce}_{\alpha}^{\rm Leib}\left(  L\right) \right)_{\rm Lie}
,\overline{\alpha}_{\rm Lie}\right)$ such that $\gamma \cdot \Psi =
u_{\alpha}$.

To end the proof we must check that $\overline{\Phi}$ and $\Psi$
are inverse to each other. Indeed,

$u_{\alpha} \cdot \left(  \overline{\Phi} \cdot \Psi\right)  =\left(  u_{\alpha} \cdot \overline{\Phi}\right) \cdot  \Psi =\gamma \cdot \Psi=u_{\alpha}$

and

$\gamma \cdot \left(  \Psi \cdot \overline{\Phi}\right)  = \left(
\gamma \cdot \Psi\right)  .\overline{\Phi}=u_{\alpha} \cdot \overline{\Phi}=\gamma$

 Since $\left(
\frak{uce}_{\alpha}^{\rm Lie}\left( L\right)
,\widetilde{\alpha} \right)$ and $\left(
\left(\frak{uce}_{\alpha}^{\rm Leib}\left(  L\right) \right)_{\rm Lie}
,\overline{\alpha}_{\rm Lie}\right)$ are perfects, then by Lemma
\ref{lema 2} in \cite{CIP}  we conclude that
$\overline{\Phi}.\Psi=\Psi.\overline{\Phi}=Id$. \rdg
\bigskip

Let $\left(  L, [-,-],\alpha_{L}\right)$ be a perfect Hom-Lie algebra.  By Theorem 4.11  $c)$ and  $d)$ in \cite{CIP}, $\left(  L, [-,-],\alpha_{L}\right)$ admits a universal central extension in the  category {\sf Hom-Lie} of the form:
$$0 {\longrightarrow} \left( H_{2}^{\alpha}\left(  L\right)
,\widetilde{\alpha}\right)  \overset{}{\longrightarrow}  \left(
\frak{uce}_{\rm Lie}\left(  L\right)  ,\widetilde{\alpha}\right)   \overset{u_{L}%
}{\longrightarrow} \left(  L,\alpha_{L}\right) {\longrightarrow} 0
$$
Here $\frak{uce}_{\rm Lie}(L)$ denotes the quotient  $\mathbb{K}$-vector space   $\frac{L \wedge L}{I_L}$, where $I_L$ is the  subspace of $L \wedge L$ spanned by the elements of the form $-[x_1,x_2] \wedge \alpha_L(x_3) + [x_1,x_3] \wedge \alpha_L(x_2) - [x_2,x_3] \wedge \alpha_L(x_1), x_1, x_2, x_3 \in L$, that is, $I_L = {\rm Im} \left( d_3 : C_3^{\alpha}(L) = \Lambda^3 L \to C_2^{\alpha}(L) = \Lambda^2 L \right)$ and $\left( H_{2}^{\alpha}\left(  L\right),\widetilde{\alpha}\right)$ denotes the second homology with trivial coefficients of  the Hom-Lie algebra $(L, \alpha_L)$ (see \cite{CIP} for details).

On the other hand, since a perfect Hom-Lie algebra is a perfect Hom-Leibniz algebra as well, then   by Theorem \ref{teorema} $c)$ and  $d)$, $\left(  L, [-,-],\alpha_{L}\right)$ also admits a universal central extension in the category {\sf Hom-Leib} of the form:
$$0 {\longrightarrow}
\left(  HL_{2}^{\alpha}\left(  L\right)
,\overline{\alpha} \right)  \overset{}{\longrightarrow}  \left(
\frak{uce}_{\rm Leib}\left(  L\right)  ,\overline{\alpha} \right)   \overset{U_{L}%
}{\longrightarrow}  \left(  L,\alpha_{L}\right)  {\longrightarrow}0
$$

\begin{Pro}
Let $\left(  L,\alpha_{L}\right)$ be a perfect Hom-Lie algebra and
let $\left( \frak{uce}_{\rm Lie}\left(  L\right), \widetilde{\alpha}\right)$ and $\left( \frak{uce}_{\rm Leib}\left( L\right),
\overline{\alpha} \right)$ be its  universal central extensions
in the categories {\sf Hom-Lie} and {\sf Hom-Leib}, respectively.
$$\left( \frak{uce}_{\rm Lie}\left(  L\right),\widetilde{\alpha} \right)  \simeq \left( \frak{uce}_{\rm Leib}\left(
L\right), \overline{\alpha} \right)  \Longleftrightarrow H_{2}^{\alpha}\left(  L\right)  \simeq
HL_{2}^{\alpha}\left(  L\right).$$
\end{Pro}
{\it Proof}. If $\left(  L,\alpha_{L}\right)$ is perfect in the
category {\sf Hom-Lie}, then it is  perfect in the category {\sf
Hom-Leib} as well. Consequently, by Theorem 4.11 {\it c)} in \cite{CIP}
and Theorem \ref{teorema} {\it c)}, $\left( L,\alpha_{L}\right)$ admits
universal central extension in both categories, hence we can
construct the following commutative diagram:
$$\xymatrix{
  0  \ar[r] & HL_2^{\alpha}(L) \ar[r] \ar@{-->}[d]^{\sigma_{\mid}} & \left(  \frak{uce}_{\rm Leib}\left(  L\right)  ,\overline{\alpha} \right) \ar[r]^{\ \ \ U} \ar@{-->}[d]^{\exists! \sigma}& \left(  L,\alpha_{L}\right) \ar[r]\ar@{=}[d] & 0 \quad ({\sf Hom-Leib})\\
 0 \ar[r] &  H_2^{\alpha}(L) \ar[r] & \left(  \frak{uce}_{\rm Lie}\left(  L\right)  ,\widetilde{\alpha} \right) \ar[r]^{\ \ \ u} & \left(  L,\alpha_{L}\right) \ar[r] & 0 \quad  ({\sf Hom-Lie})}$$

$U:\left(  \frak{uce}_{\rm Leib}\left(
L\right),\overline{\alpha} \right)  \to\left( L,\alpha_{L}\right)$
is a universal central extension in the category {\sf Hom-Leib}
and $u:\left( \frak{uce}_{\rm Lie} \left(
L\right),\widetilde{\alpha} \right) \to\left(  L,\alpha_{L}\right)$
is a central extension in same category, then there exists a
unique homomorphism $\sigma:\left( \frak{uce}_{\rm Leib}\left(
L\right),\overline{\alpha} \right)  \to \left(
\frak{uce}_{\rm Lie}\left( L\right),\widetilde{\alpha}\right)$ such
that $u \cdot \sigma=U$.

If $\sigma$ is an isomorphism, then its restriction
$\sigma_{\mid}$ is an isomorphism as well, hence
$H_{2}^{\alpha}\left( L\right)  \simeq HL_{2}^{\alpha}\left(
L\right)$.

Conversely, if $\sigma_{\mid}: HL_{2}^{\alpha}\left( L\right)  \to
H_{2}^{\alpha}\left(  L\right)$ is an isomorphism, then
 $\sigma$ is an isomorphism as well, by the the Short Five Lemma \cite{BB} and having in mind that {\sf Hom-Leib} is a semi-abelian category. \rdg

\begin{Rem}
When Leibniz (Lie) algebras are considered as Hom-Leibniz
(Hom-Lie) algebras, i.e. when $\alpha = Id$, then the above
results recover the relationship between universal central
extensions of a perfect Lie algebra in the categories of Lie
algebras and Leibniz algebras given in \cite{Gn 1}.
\end{Rem}

\bigskip \bigskip
\noindent{\bf Acknowledgements}

First and second  authors  were supported by  Ministerio de
 Ciencia e Innovación (Spain), Grant MTM2009-14464-C02 (European
FEDER support included) and by Xunta de
Galicia, Grant Incite09 207 215 PR.

\bibliographystyle{model1a-num-names}
\bibliography{<your-bib-database>}

\end{document}